\documentclass{article}
\usepackage{ascmac,amsmath,amssymb,amsthm}
\newtheorem{thm}{\bf{Theorem}}

\newtheorem{lem}{\bf{Lemma}}

\begin{document}
\title{Supremum of the function $S_1(t)$ on short intervals}
\author{Takahiro Wakasa}
\date{Graduate School of Mathematics, Nagoya University, Chikusa-ku, Nagoya 464-8602, Japan;\\ e-mail: d11003j@math.nagoya-u.ac.jp}
\maketitle
\begin{abstract}
We prove a lower bound on the supremum of the function $S_1(T)$ on short intervals, defined by the integration of the argument of the Riemann zeta-function. The same type of result on the supremum of $S(T)$ have already been obtained by Karatsuba and Korolev. Our result is based on the idea of the paper of Karatsuba and Korolev. Also, we show an improved Omega-result for a lower bound.
\end{abstract}\section{Introduction}
We consider the argument of the Riemann zeta function $\zeta(s)$, where $s=\sigma+ti$ is a complex variable, on the critical line $\sigma=\frac{1}{2}$.

We introduce the functions $S(t)$ and $S_1(t)$. When $T$ is not the ordinate of any zero of $\zeta(s)$, we define
\begin{align*}
S(T)=\frac{1}{\pi}\arg\zeta\left(\frac{1}{2}+Ti\right).
\end{align*} 
This is obtained by continuous variation along the straight lines connecting $2$, $2+Ti$, and $\frac{1}{2}+Ti$, starting with the value zero. When $T$ is the ordinate of some zero of $\zeta(s)$, we define
\begin{align*}
S(T)=\frac{1}{2}\{S(T+0)+S(T-0)\}.
\end{align*}

Next, we define $S_1(T)$ by 
\begin{align*}
S_1(T)=\int_{0}^{T}S(t)dt+C,
\end{align*}
where $C$ is the constant defined by
\begin{align*}
C=\frac{1}{\pi}\int_{\frac{1}{2}}^{\infty}\log|\zeta(\sigma)|d\sigma. 
\end{align*}

It is a classical result of von Mangoldt (cf. chapter 9 of Titchmarsh \cite{E.C.TITCHMARSH}) that there exists a number $T_0>0$ such that for $T>T_0$ we have
\begin{align*}
S(T)=O(\log T).
\end{align*}
Also, Littlewood \cite{Littlewood} proved that there exists a number $T_0>0$ such that for $T>T_0$ we have
\begin{align*}
S_1(T)=O(\log T).
\end{align*} 
Further, Littlewood proved that under the Riemann Hypothesis we have
\begin{align*}
S(T)=O\left(\frac{\log T}{\log\log T}\right)
\end{align*}
and
\begin{align*}
S_1(T)=O\left(\frac{\log T}{(\log \log T)^2}\right).
\end{align*}

 There exist some known results for $S(t)$ on short intervals. In $1946$, Selberg \cite{Selberg1} proved the inequalities
\begin{align*}
\sup_{T\leq t\leq 2T}(\pm S(t))\geq A\frac{(\log T)^\frac{1}{3}}{(\log \log T)^\frac{7}{3}},
\end{align*}
where $A$ is a positive absolute constant. Also, a similar result for $S_1(t)$ is 
\begin{align}
S_1(t)=\Omega_{\pm}\left(\frac{(\log t)^\frac{1}{3}}{(\log\log t)^\frac{10}{3}}\right).\label{1111}
\end{align}
Also, Tsang \cite{Tsang} proved for $S_1(t)$ that
\begin{align}
S_1(t)=\left\{ \begin{array}{lll}
\Omega_{+}\left(\frac{(\log t)^\frac{1}{2}}{(\log \log t)^\frac{9}{4}}\right)~~~&unconditionally,\\
\Omega_{-}\left(\frac{(\log t)^\frac{1}{3}}{(\log \log t)^\frac{4}{3}}\right)~~~&unconditionally, \\
\Omega_{\pm}\left(\frac{(\log t)^\frac{1}{2}}{(\log \log t)^\frac{3}{2}}\right)~~~&assuming~R.H. \label{111}
\end{array} \right.
\end{align}
In $1977$, Mongomery \cite{Montgomery1} established the following result under the assumption of the Riemann hypothesis: In the interval $(T^\frac{1}{6},T)$, there exist points $t_0$ and $t_1$ such that
\begin{align*}
(-1)^jS(t_j)\geq \frac{1}{20}\left(\frac{\log T}{\log\log T}\right)^\frac{1}{2}, ~~~j=0,1.
\end{align*} 
In $1986$, Tsang \cite{Tsang} improved the methods of \cite{Selberg1} to obtain the following inequalities strengthening the above results of Selberg and Mongomery:
\begin{align*}
\sup_{T\leq t \leq 2T}(\pm S(t))\geq A\left(\frac{\log T}{\log\log T}\right)^a,
\end{align*}
where $A>0$ is an absolute constant and the value of $a$ is equal to $\frac{1}{2}$ if the Riemann hypothesis is true and equal to $\frac{1}{3}$ otherwise.
 
In $2005$, Karatsuba and Korolev \cite{KandK} established the following result: Let $0<\epsilon<\frac{1}{10^3}$, $T \geq T_0(\epsilon)>0$, and $H=T^{\frac{27}{82}+\epsilon}$. Then
\begin{align*}
\sup_{T-H \leq t \leq T+2H}(\pm S(t)) \geq \frac{\epsilon^\frac{5}{4}}{1000}\left(\frac{\log T}{\log\log T}\right)^\frac{1}{3}.
\end{align*}
Our result in the present paper is obtained by applying the method of proving the above result to the function $S_1(t)$. 
\begin{thm}\label{th}
@\\~~~Let $0<\epsilon<\frac{1}{10^3}$, $T \geq T_0(\epsilon)>0$, and $H=T^{\frac{27}{82}+\epsilon}$. Then
\begin{align*}
\sup_{T-H \leq t \leq T+2H}(\pm S_1(t)) \geq \frac{\epsilon}{4000\pi}\left(\frac{(\log T)^\frac{1}{3}}{(\log\log T)^\frac{5}{3}}\right).
\end{align*}
\end{thm}
This can be proven similarly to the above result of Karatsuba and Korolev \cite{KandK}. So in this paper, we describe just the outline of the proof of Theorem $1$. However, lemmas to apply for the proof of Theorem $1$ are different from those in \cite{KandK}. There are five lemmas to apply, four lemmas among them are different. Therefore, we describe the details of the proofs of those lemmas, which are Lemma \ref{lem1}, Lemma \ref{lem2}, Lemma \ref{lem3} and Lemma \ref{lem4}. The basic ideas of the proofs of Lemmas \ref{lem1}, \ref{lem2}, \ref{lem3} and \ref{lem4} are based on the proof of Theorem $2$, Lemma $2$, Lemma $4$ and Lemma $3$, respectively, of Chapter $3$ in Karatsuba and Korolev \cite{KandK}.
\begin{thm}\label{th2}
\begin{align*}
S_1(t)=\Omega_\pm \left(\frac{(\log t)^\frac{1}{3}}{(\log\log t)^\frac{5}{3}}\right).
\end{align*}
\end{thm}
Theorem \ref{th2} can be seen immediately from Theorem \ref{th}. This is an improvement of Selberg's result (\ref{1111}). Moreover, for $\Omega_+$, Theorem \ref{th2} is also an improvement of Tsang's result (\ref{111}).

There are functions $S_2(t)$, $S_3(t),\cdots$ defined by
\begin{align*}
S_m(t)=\int_{0}^{t}S_{m-1}(u)du+C_m
\end{align*}
for $m\geq 2$, where constants $C_m$ depend on $m$. It seems that we cannot apply the method in Karatsuba and Korolev \cite{KandK} for $S_2(t)$, $S_3(t),$ etc. because $S_2(t),$ etc. do not have the expression like 
\begin{align}
S_1(t)=\frac{1}{\pi}\int_{\frac{1}{2}}^{\frac{3}{2}}\log|\zeta(\sigma+ti)|d\sigma +O(1) \label{1.0}
\end{align}
for $S_1(t)$ in p. 274 of Selberg \cite{Selberg3}. This expression is essential in the proof of Lemma \ref{lem1}. The basic idea of the method in Karatsuba and Korolev \cite{KandK} relies on Lemma \ref{lem1}. Therefore, the method in this paper cannot be applied to  $S_2(t),$ etc. 

Therefore, some new idea or the expression like (\ref{1.0}) will be necessary to obtain the result similar to our Theorem \ref{th},  for functions $S_2(t),$ etc.

\section{Some lemmas}
Here we introduce the following notations.

Let $2\leq x\leq t^2$. We set
\begin{align*}
\sigma_{x,t}=\frac{1}{2}+2\max\left(\left|\beta-\frac{1}{2}\right|,\frac{1}{\log x}\right),
\end{align*}
where $\beta$ ranges over the real parts of the zeros $\rho=\beta+\gamma i$ of the Riemann zeta function that satisfy the condition 
\begin{align*}
|\gamma-t|\leq\frac{x^{3\left|\beta-\frac{1}{2}\right|}}{\log x}.
\end{align*}
Also, we set
\begin{align*}
\Lambda(n)=\left\{ \begin{array}{ll}
\log p &~~ {\rm if}~n=p^k~{\rm with~ a~ prime}~p~{\rm and~ an~ integer}~k\geq 1, \\
0&~~{\rm otherwise}.\\
\end{array} \right.
\end{align*}
Using these notations, we state the following lemmas.
\begin{lem}\label{lem1}
　\\~~~Let $f(z)$ be a function taking real values on the real line, analytic on the strip $|\Im z| \leq 1$, and satisfying the inequality $|f(z)| \leq c(|z|+1)^{-(1+\alpha)}$, $c>0$, $\alpha >0$, on this strip. Then, the formula 
\begin{align*}
\int_{-\infty}^{\infty}f(u)S_1(t+u)du = &\frac{1}{\pi}\sum_{n=2}^{\infty}\frac{\Lambda(n)}{n^\frac{1}{2}(\log n)^2}\Re \left(\frac{1}{n^{ti}}\hat{f}(\log n)\right)-C\hat{f}(0)\\
&+2\Biggl\{\sum_{\beta > \frac{1}{2}}\int_{\frac{1}{2}}^{\beta}\int_{0}^{\beta-\sigma}\Re f(\gamma-t-xi) dx d\sigma \\
&~~~~~~~-\int_{\frac{1}{2}}^{1}\int_{0}^{1-\sigma}\Re f(-t-xi) dx d\sigma \Biggr\},
\end{align*}
where $\hat{f}(x)$ is given by the formula 
\begin{align*}
\hat{f}(x)=\int_{-\infty}^{\infty}f(u)e^{-ixu}du,
\end{align*}
holds for any $t$, where the summation in the last sum is taken over all complex zeros $\rho = \beta+ \gamma i$ of $\zeta(s)$ to the right of the critical line, and where
\begin{align*}
C=\frac{1}{\pi}\int_{\frac{1}{2}}^{\infty}\log|\zeta(\sigma)|d\sigma.
\end{align*}
\end{lem}
\begin{lem}\label{lem2}
　\\~~~For any sufficiently large positive values of $H$, $t$, and $\tau$ with $\tau < \log t$ and $H<t$, 
\begin{align*}
\int_{-\frac{1}{2}H\tau}^{\frac{1}{2}H\tau}\left(\frac{\sin u}{u}\right)^2 S_1\left(t+\frac{2u}{\tau}\right)du=W(t)+R(t)+O\left(\frac{\log t}{\tau H}\right)+O(1),
\end{align*}
where
\begin{align*}
W(t)&=\sum_{p\leq e^\tau}\frac{\cos(t\log p)}{p^\frac{1}{2}\log p}\left(1-\frac{\log p}{\tau}\right),\\
R(t)&=\tau \sum_{\beta > \frac{1}{2}}\int_{\frac{1}{2}}^{\beta}\int_{0}^{\beta-\sigma}\Re\left(\frac{\sin\frac{\tau}{2}(\gamma-t-xi)}{\frac{\tau}{2}(\gamma-t-xi)}\right)^2 dx d\sigma.
\end{align*}
\end{lem}
\begin{lem}\label{lem3}
　\\~~~Let $\epsilon$ with $0<\epsilon<\frac{1}{1000}$ be fixed. Let $T\geq T_0(\epsilon)>0$, $H=T^{\frac{27}{82}+\epsilon}$ and $k$ be an integer such that $k\geq k_0(\epsilon)>1$, let $m=2k+1$, $\tau=2\log\log H$, and $m\tau<\frac{1}{10}\epsilon \log T$. Then the function $R(t)$ defined by Lemma \ref{lem2} satisfies the inequality
\begin{align*}
\int_{T}^{T+H}|R(t)|^mdt<H\left\{25^m+(\log T)^3\left(\frac{50\tau m^2}{\epsilon ^3\log T}\right)^m\right\}.
\end{align*}
\end{lem}

\begin{lem}\label{lem4}
　\\~~~Let $T\geq T_0 >0$, $e^2<H<T$, $2<\tau <\log H$, and $k$ be an integer such that $k\geq k_0 >1$ and $(2k\log k)^2<e^{\frac{4}{5}\tau}$. Then
\begin{align}
&\int_{T}^{T+H}W(t)^{2k}dt>\left(\frac{1}{5\sqrt{10e}}\cdot\frac{k^\frac{1}{2}}{\log k}\right)^{2k}H-e^{3k\tau},\label{lem4.1}\\
&\left|\int_{T}^{T+H}W(t)^{2k+1}dt\right|<e^{3k\tau+\frac{3}{2}\tau}.\label{lem4.2}
\end{align}
\end{lem}
This lemma is Lemma $3$ of Chapter $3$ in Karatsuba and Korolev \cite{KandK}. But in Karatsuba and Korolev \cite{KandK}, the function $W(t)$ is defined by 
\begin{align*}
W(t)=-\sum_{p\leq e^\tau}\frac{\sin(t\log p)}{p^\frac{1}{2}}\left(1-\frac{\log p}{\tau}\right),
\end{align*}
which are defferent from the definition in this paper.

The following lemma is given in Karatsuba and Korolev \cite{KandK}.
\begin{lem}\label{lem5}
　\\~~~Let $H>0$ and $M>0$, let $k\geq 1$ be an integer, and let $W(t)$, $R(t)$ be real functions which satisfy the conditions
\begin{align*}
&1)~~\int_{T}^{T+H}|W(t)|^{2k}dt\geq HM^{2k},\\
&2)~~\left|\int_{T}^{T+H}W(t)^{2k+1}dt\right|\leq \frac{1}{2}HM^{2k+1},\\
&3)~~\int_{T}^{T+H}|R(t)|^{2k+1}dt<H\left(\frac{M}{2}\right)^{2k+1}.
\end{align*}
Then
\begin{align*}
\max_{T\leq t\leq T+H}\pm(W(t)+R(t))\geq \frac{1}{8}M.
\end{align*}
\end{lem}
This lemma is Lemma $1$ of Chapter $3$ in Karatsuba and Korolev \cite{KandK}.
\section{Proof of Lemma \ref{lem1}}
This proof is an analogue of the proof of Theorem $2$ of Chapter $3$ in Karatsuba and Korolev \cite{KandK}.

\textit{Proof}.
Put $\frac{1}{2}\leq \sigma \leq\frac{3}{2}$. We set $\psi(z)=f((\sigma-z)i-t)$ and take $X\geq 2(|t|+10)$ such that the distance from the ordinate of any zero of $\zeta(s)$ to $X$ is not less than $c(\log X)^{-1}$, where $c$ is a positive absolute constant.

Let $\Gamma$ be the boundary of the rectangle with the vertices $\sigma \pm Xi$, $\frac{3}{2}\pm Xi$, and let a horizontal cut be drawn from the line $\Re s=\sigma$ inside this rectangle to each zero $\rho=\beta+\gamma i$ and also to the point $z=1$. Then the functions $\log \zeta(z)$ and $\psi(z)$ are analytic inside $\Gamma$. 

By the residue theorem, the following equality holds:
\begin{align*}
0&=\int_{\Gamma}\psi(z)\log \zeta(z)dz \\
 &=\left(\int_{\frac{3}{2}-Xi}^{\frac{3}{2}+Xi}-\int_{\sigma+Xi}^{\frac{3}{2}+Xi}-\int_{\sigma-Xi}^{\sigma+Xi}+\int_{\sigma-Xi}^{\frac{3}{2}-Xi}\right) \psi(z) \log \zeta(z)dz \\
 &=I_1-I_2-I_3+I_4,
\end{align*}
say. Then, we have
\begin{align*}
I_1=i\sum_{n=2}^{\infty}\frac{\Lambda(n)}{n^{\sigma+ti}\log n}\hat{f}(\log n)+O\left(\frac{1}{X^\alpha}\right)
\end{align*}
since for $\alpha =\frac{3}{2}-\sigma$
\begin{align*}
\int_{-\infty}^{\infty}\psi\left(\frac{3}{2}+ui\right)\log\zeta\left(\frac{3}{2}+ui\right)du&=\sum_{n=2}^{\infty}\frac{\Lambda(n)}{n^\frac{3}{2}\log n}\int_{-\infty}^{\infty}\frac{1}{n^{ui}}f(u-t-\alpha i)du\\
&=\sum_{n=2}^{\infty}\frac{\Lambda(n)}{n^{\sigma+ti}\log n}\hat{f}(\log n).
\end{align*}
Also,
\begin{align*}
I_2&=O\left(\frac{(\log X)^2}{X^{(1+\alpha)}}\right),\\
I_4&=O\left(\frac{(\log X)^2}{X^{(1+\alpha)}}\right)
\end{align*}
as in p. $461$ of Karatsuba and Korolev \cite{KandK}.

We denote by $L$ the cut going from the point $\sigma+\gamma i$ to the each points $\beta+\gamma i$, and denote by $I(L)$ the integral over the banks of this cut. Then,  
\begin{align*}
I(L)=\int_{L}\psi(z)\log\zeta(z)dz=2\pi i\sum_{ \beta > \sigma }\int_{0}^{\beta-\sigma}f(\gamma-t-xi)dx
\end{align*}
as in p. $462$ of Karatsuba and Korolev \cite{KandK}.

If $L$ is the cut going to the point $z=1$, then
\begin{align*}
I(L)=-2\pi i \int_{0}^{1-\sigma}f(-t-xi)dx.
\end{align*}  
Hence, we have 
\begin{align*}
I_3&=\int_{\sigma-Xi}^{\sigma+Xi}\psi(z)\log\zeta(z)dz\\
&~~~~~~~~~-2\pi i\left(\sum_{\beta > \sigma }\int_{0}^{\beta-\sigma}f(\gamma-t-xi)dx-\int_{0}^{1-\sigma}f(-t-xi)dx\right).
\end{align*}
When $X$ tends to infinity, we obtain
\begin{align*}
\lim_{X\rightarrow \infty}\int_{\sigma-Xi}^{\sigma+Xi}\psi(z)\log\zeta(z)dz &=i\int_{-\infty}^{\infty}f(u)\log\zeta(\sigma+(t+u)i)du\\
&=i\sum_{n=2}^{\infty}\frac{\Lambda(n)}{n^{\sigma+ti}\log n}\cdot \hat f(\log n)\\
&~~+2\pi i\left(\sum_{ \beta > \sigma}\int_{0}^{\beta-\sigma}f(\gamma-t-xi)dx-\int_{0}^{1-\sigma}f(-t-xi)dx\right).
\end{align*}
Dividing by $i$, we get for $\sigma \geq \frac{1}{2}$
\begin{align*}
\int_{-\infty}^{\infty}f(u)& \log\zeta(\sigma+(t+u)i)du =\sum_{n=2}^{\infty}\frac{\Lambda(n)}{n^{\sigma+ti}\log n}\hat f(\log n)\\
&+2\pi \left(\sum_{\beta > \sigma }\int_{0}^{\beta-\sigma}f(\gamma-t-xi)dx-K(\sigma)\int_{0}^{1-\sigma}f(-t-xi)dx\right),
\end{align*}
where
\begin{align*}
K(\sigma)=\left\{ \begin{array}{ll}
1&~~ {\rm for}~\frac{1}{2}\leq \sigma\leq 1, \\
0&~~{\rm for}~\sigma>1.\\
\end{array} \right.
\end{align*}

Here, taking the real part and applying (\ref{1.0}) and integrating in $\sigma$ over the interval $[\frac{1}{2},\frac{3}{2}]$, we have
\begin{align*}
&\int_{-\infty}^{\infty}\int_{\frac{1}{2}}^{\frac{3}{2}}f(u) \log|\zeta(\sigma+(t+u)i)|d\sigma du\\
&=\pi\int_{-\infty}^{\infty}S_1(t+u)f(u)du+\pi\int_{-\infty}^{\infty}f(u) Cdu\\
&=\sum_{n=2}^{\infty}\frac{\Lambda(n)}{n^{\frac{1}{2}}(\log n)^2}\Re\left(\frac{1}{n^{ti}}\hat f(\log n)\right)\\
&~~~+2\pi\left(\int_{\frac{1}{2}}^{\frac{3}{2}}\sum_{ \beta > \sigma}\int_{0}^{\beta-\sigma} \Re f(\gamma-t-xi)dxd\sigma-\int_{\frac{1}{2}}^{\frac{3}{2}}\int_{0}^{1-\sigma}\Re f(-t-xi)dxd\sigma\right).
\end{align*}
Therefore
\begin{align*}
&\int_{-\infty}^{\infty}S_1(t+u)f(u)du\\
&=\frac{1}{\pi}\sum_{n=2}^{\infty}\frac{\Lambda(n)}{n^{\frac{1}{2}}(\log n)^2}\Re\left(\frac{1}{n^{ti}}\hat f(\log n)\right)-C\hat{f}(0)\\
&~~+2\left(\int_{\frac{1}{2}}^{\frac{3}{2}}\sum_{ \beta > \sigma }\int_{0}^{\beta-\sigma} \Re f(\gamma-t-xi)dxd\sigma-\int_{\frac{1}{2}}^{1}\int_{0}^{1-\sigma}\Re f(-t-xi)dxd\sigma\right)\\
&=\frac{1}{\pi}\sum_{n=2}^{\infty}\frac{\Lambda(n)}{n^{\frac{1}{2}}(\log n)^2}\Re\left(\frac{1}{n^{ti}}\hat f(\log n)\right)-C\hat{f}(0)\\
&~~+2\left(\sum_{ \beta > \frac{1}{2}}\int_{\frac{1}{2}}^{\beta}\int_{0}^{\beta-\sigma}\Re f(\gamma-t-xi)dxd\sigma-\int_{\frac{1}{2}}^{1}\int_{0}^{1-\sigma}\Re f(-t-xi)dxd\sigma\right).
\end{align*}
\qed
\section{Proof of Lemma \ref{lem2}}
This proof is an analogue of the proof of Lemma $2$ of Chapter $3$ in Karatsuba and Korolev \cite{KandK}.

\textit{Proof}.
Put $f(z)=\left(\frac{\sin\frac{\tau z}{2}}{\frac{\tau z}{2}}\right)^2$. By
\begin{align*}
\hat{f}(x)=\int_{-\infty}^{\infty}e^{-xui}f(u)du=\frac{2\pi}{\tau}\max\left(0,1-\left|\frac{x}{\tau}\right|\right),
\end{align*}
we get
\begin{align*}
\hat{f}(\log n)=\left\{\begin{array}{ll}
\frac{2\pi}{\tau}\left(1-\frac{\log n}{\tau}\right) & (1\leq n\leq e^\tau)\\
0 & (n>e^\tau)
\end{array}\right. .
\end{align*}
Then, we have 
\begin{align}
\int_{-\infty}^{\infty}\left(\frac{\sin \frac{\tau u}{2}}{\frac{\tau u}{2}}\right)^2 S_1(t+u)du
&=\frac{1}{\pi}\sum_{n\leq e^\tau}\frac{\Lambda(n)}{n^{\frac{1}{2}}(\log n)^2}\cdot \frac{2\pi}{\tau}\left(1-\frac{\log n}{\tau}\right)\cos(t \log n)\nonumber \\
&~~~+2\Biggl\{\sum_{\beta > \frac{1}{2} }\int_{\frac{1}{2}}^{\beta}\int_{0}^{\beta-\sigma}\Re\left(\frac{\sin \frac{\tau}{2}(\gamma-t-\xi i)}{\frac{\tau}{2}(\gamma-t-\xi i)}\right)^2d\xi d\sigma \nonumber \\
&~~~~~~~~-\int_{\frac{1}{2}}^{1}\int_{0}^{1-\sigma}\Re\left(\frac{\sin \frac{\tau}{2}(\gamma-t-\xi i)}{\frac{\tau}{2}(\gamma-t-\xi i)}\right)^2 d\xi d\sigma \Biggr\}\nonumber \\
&~~~-C\hat{f}(0)\label{lem2-1}
\end{align}
by Lemma \ref{lem1}. Since for $0\leq \xi \leq 1-\sigma$
\begin{align*}
2\left|\frac{\sin \frac{\tau}{2}(t+\xi i)}{\frac{\tau}{2}(t+\xi i)}\right|^2<\frac{1}{5\tau}
\end{align*}
as in p. $473$ of Karatsuba and Korolev \cite{KandK}, we have 
\begin{align*}
\left|2\int_{\frac{1}{2}}^{1}\int_{0}^{1-\sigma}\Re\left(\frac{\sin \frac{\tau}{2}(\gamma-t-\xi i)}{\frac{\tau}{2}(\gamma-t-\xi i)}\right)^2 d\xi d\sigma \right|=O\left(\frac{1}{\tau}\right).
\end{align*}
In the first term of the right-hand side in (\ref{lem2-1}), we single out the terms corresponding to the $n=p$  in the sum and estimate the remainder terms. Then, we have 
\begin{align*}
\sum_{2\leq k} \sum_{p^k \leq e^\tau} \frac{\Lambda(p^\frac{k}{2})\cos(t\log p^k)}{p^\frac{k}{2}(\log p^k)^2}\cdot\frac{2}{\tau}\left(1-\frac{\log p^k}{\tau}\right)<\sum_{2\leq k} \sum_{p^k \leq e^\tau}\frac{\log p}{p^\frac{k}{2}(\log p^k)^2}\cdot\frac{2}{\tau}\ll\frac{1}{\tau}.
\end{align*} 
Hence
\begin{align}
\int_{-\infty}^{\infty}\left(\frac{\sin\frac{\tau u}{2}}{\frac{\tau u}{2}}\right)^2 S_1(t+u)du&=\frac{2}{\tau}\sum_{p\leq e^\tau}\frac{\cos(t\log p)}{p^\frac{1}{2}\log p}\left(1-\frac{\log p}{\tau}\right)-C\cdot\frac{2\pi}{\tau}\nonumber \\
&~~+2\sum_{ \beta > \frac{1}{2} }\int_{\frac{1}{2}}^{\beta}\int_{0}^{\beta-\sigma}\Re\left(\frac{\sin\frac{\tau}{2}(\gamma-t-\xi i)}{\frac{\tau}{2}(\gamma-t-\xi i)}\right)d\xi d\sigma \nonumber \\
&~~+O\left(\frac{1}{\tau}\right). \label{lem2-2}
\end{align}
Put $v=\frac{\tau u}{2}$. Then the left-hand side of the above is equal to
\begin{align*}
\left(\int_{-\frac{1}{2}H\tau}^{\frac{1}{2}H\tau}+\int_{-\infty}^{-\frac{1}{2}H\tau}+\int_{\frac{1}{2}H\tau}^{\infty}\right)\left(\frac{\sin v}{v}\right)^2S_1\left(t+\frac{2v}{\tau}\right)\frac{2}{\tau}dv.
\end{align*}
Since $S_1(t)=O(\log t)$, we have
\begin{align*}
\left|\left(\int_{-\infty}^{-\frac{1}{2}H\tau}+\int_{\frac{1}{2}H\tau}^{\infty}\right)\left(\frac{\sin v}{v}\right)^2S_1\left(t+\frac{2v}{\tau}\right)dv\right| &\ll\frac{1}{\tau}\int_{H}^{\infty}\log(t+v')\frac{1}{v'^2}dv'\\
&\ll \frac{1}{\tau}\left\{\int_{H}^{t}\frac{\log t}{v'^2}dv'+\int_{t}^{\infty}\frac{\log v'}{v'^2}dv'\right\}\\
&\ll \frac{1}{\tau}\left(\frac{\log t}{H}+\frac{\log t}{t}\right)\ll\frac{\log t}{\tau H}.
\end{align*}
Inserting these estimates into (\ref{lem2-2}) and dividing by $\frac{2}{\tau}$ the both sides, we obtain the result.

\qed
\section{Proof of Lemma \ref{lem3}}
This proof is an analogue of the proof of Lemma $4$ of Chapter $3$ in Karatsuba and Korolev \cite{KandK}.

\textit{Proof}.
We put
\begin{align*}
L_k=\int_{T}^{T+H}|R(t)|^{2k+1}dt
\end{align*}
and note the inequality
\begin{align*}
\left|\Re\left(\frac{\sin(x-yi)}{x-yi}\right)^2\right|<\frac{8ye^{2y}}{1+x^2+y^2}
\end{align*}
for any $x$, $y\in \mathbb{R}$, $y\geq 0$ similarly to pp. $476-477$ of Karatsuba and Korolev \cite{KandK}. Then,
\begin{align*}
|R(t)|&\leq \tau\left|\sum_{\begin{subarray}{c} \gamma \\ \beta > \frac{1}{2} \end{subarray}}\int_{\frac{1}{2}}^{\beta}\int_{0}^{\beta-\sigma}\Re\left(\frac{\sin\frac{\tau}{2}(\gamma-t-\xi i)}{\frac{\tau}{2}(\gamma-t-\xi i)}\right)^2d\xi d\sigma \right|\\
&\leq \tau \sum_{\begin{subarray}{c} \gamma \\ \beta > \frac{1}{2} \end{subarray}}\int_{\frac{1}{2}}^{\beta}\int_{0}^{\beta-\sigma}\frac{8\cdot \frac{\tau \xi}{2}e^{\tau \xi}}{1+\left\{\frac{\tau}{2}(\gamma-t)\right\}^2+\left(\frac{\tau \xi}{2}\right)^2}d\xi d\sigma\\
&<4\tau^2\sum_{\begin{subarray}{c} \gamma \\ \beta > \frac{1}{2} \end{subarray}}\int_{\frac{1}{2}}^{\beta}\int_{0}^{\beta-\frac{1}{2}}\frac{\xi e^{\tau(\beta-\frac{1}{2})}}{1+\left\{\frac{\tau}{2}(\gamma-t)\right\}^2+\left(\frac{\tau}{2}\left(\beta-\frac{1}{2}\right)\right)^2}d\xi d\sigma\\
&=8\sum_{\begin{subarray}{c} \gamma \\ \beta > \frac{1}{2} \end{subarray}}\left(\beta-\frac{1}{2}\right)^3\frac{e^{\tau \left(\beta-\frac{1}{2}\right)}}{\left(\frac{2}{\tau}\right)^2+(\gamma-t)^2+\left(\beta-\frac{1}{2}\right)^2}.
\end{align*}
We split the last sum into two sums. The first sum $\Sigma_1$ is the sum of the terms satisfying $|\gamma-t|>(\log T)^2$, and the second sum $\Sigma_2$ is the sum of the other terms.

Here, we denote by $\theta_t$ the largest difference of the form $\beta-\frac{1}{2}$ for zeros $\rho=\beta+\gamma i$ in the rectangle $\frac{1}{2}< \beta \leq 1$, $|\gamma-t|\leq (\log T)^2$. Also, we denote by $\theta'_t$ the supremum of the form $\beta-\frac{1}{2}$ for zeros $\rho=\beta+\gamma i$ in the rectangle $\frac{1}{2}< \beta \leq 1$, $|\gamma-t|> (\log T)^2$.

As in p. $478$ of Karatsuba and Korolev \cite{KandK}, we apply the estimation related to $\sigma_{x,t}$ and the result $N(t+1)-N(t)<18\log t$ which is obtained by the Riemann-von Mangoldt formula and $|S(t)|<8\log t$ for $t\geq t_0>0$. Then we take $x=(\log T)^\frac{1}{2} $, and we have 
\begin{align*}
\Sigma_1&<\left(\beta-\frac{1}{2}\right)\sum_{|\gamma-t|>(\log T)^2}\frac{2e^\frac{\tau}{2}}{(\gamma-t)^2}<\frac{2}{3}\theta_t'\log T\sum_{|\gamma-t|>(\log T)^2}\frac{1}{n^2}\sum_{n<|\gamma-t|\leq n+1}1\\
&<\frac{2}{3}\theta_t'\log T\cdot 36\sum_{|\gamma-t|>(\log T)^2}\frac{\log T +\log n}{n^2}<25\theta_t'
\end{align*}
and
\begin{align*}
\Sigma_2&<8\theta^3e^{\tau \theta}\sum_{|\gamma-t|\leq(\log T)^2}\frac{1}{\left(\frac{2}{\tau}\right)^2+(\gamma-t)^2+\left(\beta-\frac{1}{2}\right)^2}\\
&<8\theta^3e^{\tau \theta}\sum_{\rho}\frac{1}{(\sigma_{x,t}-\beta)^2+(\gamma-t)^2}
<8\theta^3e^{\tau \theta}\frac{13}{5}\cdot\frac{1}{\sigma_{x,t}-\frac{1}{2}}\log T\\
&\leq 8\theta^3e^{\tau \theta}\frac{13}{5}\cdot\frac{5\tau}{39}\log T=8\theta^3e^{\tau \theta}\cdot\frac{\tau}{3}\log T.
\end{align*}
From the definitions of $\theta_t$ and $\theta_t'$, we get $\theta_t<\frac{1}{2}$ and $\theta_t'<\frac{1}{2}$. Hence, we have 
\begin{align*}
|R(t)|<25\left(\theta_t'+\frac{7}{2}\theta_t^3e^{\tau\theta_t}\tau\log T\right)
<\frac{25}{2}\left(1+\frac{7}{2}\theta_t^2e^{\tau\theta_t}\tau\log T\right).
\end{align*}
Hence
\begin{align*}
L_k<\left(\frac{25}{2}\right)^m\int_{T}^{T+H}\left(1+\frac{7}{2}\theta_t^2e^{\tau\theta_t}\tau\log T\right)^mdt.
\end{align*}
This integrand is the same as that in p.~$479$ of Karatsuba and Korolev \cite{KandK}. Hence the estimation of the last integral is the same as in pp. $480-481$ of Karatsuba and Korolev \cite{KandK}. Along that way, we have
\begin{align*}
L_k &< 25^m H \left\{1+\frac{24}{5} \cdot \frac{1}{m} (\log T)^3 (2m)! \left(\frac{7}{2}\tau \log T \right)^m \left(\frac{\epsilon}{10} \log T \right)^{-2m}\right\}\\
&<25^m H\left\{1+(\log T)^3\left(\frac{2m^2\tau}{\epsilon^3 \log T}\right)^m \right\}\\
&<H\left(25^m+ (\log T)^3\left(\frac{50m^2\tau}{\epsilon^3 \log T}\right)^m \right).
\end{align*}
\qed
\section{Proof of Lemma \ref{lem4}}
This proof is an analogue of the proof of Lemma $3$ of Chapter $3$ in Karatsuba and Korolev \cite{KandK}.

\textit{Proof}.
As in pp. $474-475$ of Karatsuba and Korolev \cite{KandK}, we can write
\begin{align*}
\int_{T}^{T+H}W(t)^{2k}dt=I_k=\begin{pmatrix}
2k\\
k\\
\end{pmatrix}
\frac{H}{2^{2k}}\Sigma + \theta e^{3k \tau},
\end{align*}
where
\begin{align*}
\Sigma=\sum_{\begin{subarray}{c} p_1\cdots p_k=q_1\cdots q_k \\ p_1,\cdots, q_k\leq e^\tau \end{subarray}}f(p_1)^2\cdots f(p_k)^2,~~~f(p)=\frac{1}{p^\frac{1}{2}\log p}\left(1-\frac{\log p}{\tau}\right).
\end{align*}
Then, 
\begin{align*}
\Sigma&\geq k! \sum_{\begin{subarray}{c} p_1,\cdots, p_k~{\rm are}~{\rm distinct} \\ p_1,\cdots, p_k\leq e^\tau \end{subarray}}f(p_1)^2\cdots f(p_k)^2\\
&\geq k!\sum_{p_1\leq e^\tau}f(p_1)^2 \sum_{\begin{subarray}{c} p_2\leq e^\tau \\ p_1\neq p_2 \end{subarray}}f(p_2)^2\cdots\sum_{\begin{subarray}{c} p_k\leq e^\tau \\ p_1,\cdots, p_{k-1}\neq p_k \end{subarray}}f(p_k)^2.
\end{align*}
Since $\frac{d}{dp} f(p)^2<0$, $f(p)^2$ is monotonically decreasing function for $p\geq 2$. Also, since $(k-1)$th prime does not exceed $2k\log k$, the inner sum of the above inequality is greater than the same sum over $2k\log k<p_k<e^{\frac{4}{5}\tau}$. Hence the inner sum over $p_k$ is greater than 
\begin{align*}
\left(\frac{1}{5}\right)^2 \sum_{2k \log k<p \leq e^{\frac{4}{5}\tau}}\frac{1}{p(\log p)^2}.
\end{align*}
For $(2k\log k)^2\leq e^{\frac{4}{5}\tau}$, since 
\begin{align*}
\sum_{U<p\leq U^2}\frac{1}{p(\log p)^2}&\geq \frac{1}{4(\log U)^2}\sum_{U<p\leq U^2}\frac{1}{p}\\
&=\frac{1}{4(\log U)^2}(\log\log U^2-\log\log U+o(1))>\frac{1}{8(\log U)^2},
\end{align*}
the sum over $p_k$ is greater than $\frac{1}{10}\left(\frac{1}{5}\right)^2\frac{1}{(\log k)^2}$. Also, the same lower bound holds for the sums over $p_1,p_2,\cdots,p_{k-1}$. Therefore, we see
\begin{align*}
\Sigma \geq k!\left(\frac{1}{250(\log k)^2}\right)^k \geq \sqrt{2\pi k}\left(\frac{1}{5\sqrt{10e}}\cdot\frac{k^\frac{1}{2}}{\log k}\right)^{2k}.
\end{align*}
So, 
\begin{align*}
I_k> H\left(\frac{1}{5\sqrt{10e}}\cdot\frac{k^\frac{1}{2}}{\log k}\right)^{2k}-e^{3k\tau}.
\end{align*}
This is the first part of Lemma \ref{lem4}. The second part is proved similarly to \cite{KandK} .

\qed
\section{Outline of the proof of the Theorem \ref{th}}
As described in section $1$, our result can be proven similarly to Theorem $5$ of Chapter $3$ in Karatsuba and Korolev \cite{KandK}. Therefore, we describe the outline of the proof.

\textit{Outline of the proof}. Put $\tau=2\log\log H$. Consider the right-hand side of the inequality in the statement of Lemma \ref{lem3}. We see easily that
\begin{align*}
\frac{50\tau m^2}{\epsilon ^3\log T}<\frac{500k^2}{\epsilon^3}\cdot\frac{\log \log T}{\log T}\leq \frac{k^\frac{1}{2}}{\log k}\cdot\frac{500k^\frac{3}{2}}{\epsilon^3 }\cdot\frac{(\log\log T)^2}{\log T}=\frac{k^\frac{1}{2}}{\log k}\delta,
\end{align*}
say.

Here, putting
$
k=\left[\frac{\epsilon^2}{1000}\left(\frac{(\log T)^\frac{2}{3}}{(\log\log T)^\frac{4}{3}}\right)\right],
$
we have $\delta<\frac{1}{60}$, $(2k\log k)^2<e^{\frac{4}{5}\tau}$ and $e^{3k\tau}<H^\frac{1}{2}$. Hence, we can apply Lemma \ref{lem3} and Lemma \ref{lem4}. Then we have
\begin{align*}
&\int_{T}^{T+H}W(t)^{2k}dt>HM^{2k},\\
&\left|\int_{T}^{T+H}W(t)^{2k+1}dt\right|<\frac{1}{2}HM^{2k+1},\\
&\int_{T}^{T+H}|R(t)|^{2k+1}dt<H\left(\frac{M}{2}\right)^{2k+1},
\end{align*}
with $M=\frac{k^\frac{1}{2}}{30\log k}$. Thus, we see that $W(t)$ and $R(t)$ satisfy the conditions of Lemma \ref{lem5} with $M=\frac{k^\frac{1}{2}}{30\log k}$. Hence there are two points $t_0$ and $t_1$ such that 
\begin{align*}
W(t_0)+R(t_0)\geq \frac{M}{8}, ~~~W(t_1)+R(t_1)\leq -\frac{M}{8}
\end{align*}
in the interval $T\leq t\leq T+H$. By Lemma \ref{lem2}, we have 
\begin{align*}
&\int_{-\frac{1}{2}H\tau}^{\frac{1}{2}H\tau}\left(\frac{\sin u}{u}\right)^2 S_1\left(t_0+\frac{2u}{\tau}\right)du\geq \frac{M}{8}+O\left(\frac{\log t_0}{\tau H}\right),\\
&\int_{-\frac{1}{2}H\tau}^{\frac{1}{2}H\tau}\left(\frac{\sin u}{u}\right)^2 S_1\left(t_1+\frac{2u}{\tau}\right)du\leq -\frac{M}{8}+O\left(\frac{\log t_1}{\tau H}\right).
\end{align*}
Here, putting
\begin{align*}
M_0=\sup_{T-H\leq t\leq T+2T}S_1(t),~~~M_1=\inf_{T-H\leq t\leq T+2T}S_1(t),
\end{align*}
we have
\begin{align*}
&\int_{-\frac{1}{2}H\tau}^{\frac{1}{2}H\tau}\left(\frac{\sin u}{u}\right)^2 S_1\left(t_0+\frac{2u}{\tau}\right)du< M_0\int_{-\infty}^{\infty}\left(\frac{\sin u}{u}\right)^2=\frac{\pi}{2}M_0~~~(M_0>0),\\
&\int_{-\frac{1}{2}H\tau}^{\frac{1}{2}H\tau}\left(\frac{\sin u}{u}\right)^2 S_1\left(t_1+\frac{2u}{\tau}\right)du>M_1\int_{-\infty}^{\infty}\left(\frac{\sin u}{u}\right)^2=\frac{\pi}{2}M_1~~~(M_1<0).
\end{align*}
Therefore, we obtain for $r=0,1$
\begin{align*}
(-1)^rM_r>\frac{2}{\pi}\cdot\frac{M}{8}+O\left(\frac{\log t_r}{\tau H}\right)>\frac{1}{4\pi}\cdot\frac{k^\frac{1}{2}}{30\log k}>\frac{\epsilon}{4000\pi}\left(\frac{(\log T)^\frac{1}{3}}{(\log \log T)^\frac{5}{3}}\right).
\end{align*}
Thus, we obtain the result.\\
\\
{\bf Acknowledgments}\\

I thank my advisor Prof. Kohji Matsumoto for his advice and patience during the preparation of this paper. I also thank members in the same study, who gave adequate answers to my questions. Finally, I thank the referee who indicates errors in this paper.

\end{document}